\documentclass[11pt]{article}
\usepackage{amssymb,amsthm,amsmath}
\usepackage{latexsym}
\usepackage{graphicx}

\newtheorem{theorem}{Theorem}[section]

\newtheorem{thm}{Theorem}

\newtheorem{obs}[theorem]{Observation}

\title{Separating families of convex sets}

\date{}

\author{{D\'aniel Gerbner}\thanks{Supported by Hungarian Science Foundation EuroGIGA Grant OTKA NN 102029} \ and
{G\'eza T\'oth}\thanks{           Supported by Hungarian Science Foundation EuroGIGA Grant OTKA NN 102029
and by OTKA K 83767.} \\
R\'enyi Institute, Hungarian Academy of Sciences\\}

\begin{document}

\maketitle

\begin{abstract}
Two elements, $x$ and $y$, are separated by a set $S$
if it contains exactly one of $x$ and $y$.
We prove that any set of $n$ points in general position in the plane
can be separated by $O(n\log\log n/\log n)$ convex sets,
and for some point sets $\Omega(n/\log n)$ convex sets are necessary.
\end{abstract}

\section{Introduction}

We say that a set $S$ {\em separates} elements $x$ and $y$,
if exactly one of $x$ and $y$ is in $S$.
Given an underlying set $X$, a family $\mathcal{F}$ of its subsets
is called \emph{separating}, if for any
$x,y\in X$, $x\neq y$,  there exists an $F\in \mathcal{F}$ which
separates them.

Separating families are important tools in search theory. Suppose
there is an unknown defective element in $X$, and we can test
subsets of $X$ if they contain the
defective element or not.
We want to choose a family of
sets in advance (non-adaptively), such that testing all of its
members determines the defective element. It is not hard to see that
the family satisfies this property if and only if it is
separating. The usual goal is to test as few sets as possible,
i.e. find a separating family of minimum cardinality.

It is clear
that a separating family of $X$, $|X|=n$, contains at least
$\lceil \log n\rceil$ sets,
since $k$ subsets of $X$ divide it into at most $2^k$ parts.
(All logarithms in this paper are of base 2.)

On the other hand, we can
represent the elements of $X$ by $0-1$ sequences of length $k=\lceil \log
n\rceil$. Then, for $1\le k\le m$,
let $A_i$ denote the set of elements with 1 at coordinate $i$. Then the sets
$A_i$ form a
separating family.

Another well-known observation is the following, originally due to Bondy \cite{B72},
see also \cite{W09}.

\begin{obs}\label{mini} Let $\mathcal{F}$ be a
minimal separating family (in the sense that no proper subfamily of it is
separating). Then $|\mathcal{F}| \le n-1$.
\end{obs}

There are several different versions of this simple concept, for a survey
see \cite{DH94}. An obvious idea is to consider adaptive
algorithms, where we test only one set at a time, and choose the
next set knowing the result of the previous test. Or one can have
more defective elements, different types of tests, or errors in the
results of the tests.
One of the most studied generalizations is
that not every subset can be asked. Instead, we are given a family
$\mathcal{A}$ of subsets of $X$, and we can only test its members. A well-known
example for this is to find the defective coins, or to
sort a finite set, using only pairwise
comparisons.

In this note, our underlying set $X$ is a set of points in the plane,
and we have certain geometric restrictions
on the subsets we can ask.



For a family $\mathcal{A}$ of planar sets and a point set $X$ let
$\mathcal{A}_X=\{A \cap X | A \in \mathcal{A}\}$. For simplicity, we will call
$\mathcal{A}$ separating (with respect to $X$) if $\mathcal{A}_X$ is
separating.

\bigskip

\noindent {\bf Definition.} Let $X$ be a set of $n$
points in the plane, and let $\mathcal{A}$ be a family of planar
sets.
Let $s(X,\mathcal{A})$ denote the size of the smallest subfamily
$\mathcal{A}'\subseteq \mathcal {A}$ with the property that
the family $\{ A\cap X | A\in \mathcal{A}' \}$ is separating.
If there is no such subfamily, then let $s(X,\mathcal{A})=\infty$.

Let $s(n,\mathcal{A})$ be the maximum of $s(X,\mathcal{A})$
over all $n$--element point
sets $X$, and let $s'(n,\mathcal{A})$ be the maximum of $s(X,\mathcal{A})$
over all $n$--element point
sets $X$ in {\em general position} (that is, not three of its points are on a
line).

\bigskip

By Observation \ref{mini}, for any  family $\mathcal{A}$,
$s(n,\mathcal{A})\le n-1$ or $s(n,\mathcal{A})=\infty$, and similarly,
$s'(n,\mathcal{A})\le n-1$ or $s'(n,\mathcal{A})=\infty$.

For most of the natural families of planar
sets  $\mathcal{A}$,
it is not hard (or sometimes trivial) to give a linear lower bound for
both $s(n,\mathcal{A})$ and
$s'(n,\mathcal{A})$, and in many cases we can determine their exact values.
We only give two examples here.

\begin{thm}\label{felsikkor} Let $\mathcal{H}$ and $\mathcal{D}$
denote the family of the halfplanes and discs, respectively.
Then we have

\noindent (i) $s(n,\mathcal{H})=n-1$, $s'(n,\mathcal{H})={\lceil n/2\rceil}$,

\smallskip

\noindent (ii) $s(n,\mathcal{D})=s'(n,\mathcal{D})={\lceil n/2\rceil}$.

\end{thm}

%

%
%

The case when $\mathcal{A}$ is
the family of {\em convex} sets seems to be the most interesting.

\begin{thm}\label{fo} Let $\mathcal{A}$ denote the family of planar convex sets.
Then we have

\noindent (i) $s(n,\mathcal{A})={\lceil n/2\rceil}$,

\smallskip

\noindent and (ii) $n/2\log n\le s'(n,\mathcal{A})\le 20n\log\log n/\log n$.

\end{thm}

We prove Theorem \ref{felsikkor} in Section \ref{2}, and Theorem \ref{fo} in Section \ref{3}.

\section{Some simple families of planar sets; Proof of Theorem \ref{felsikkor}}\label{2}

It is easy to see that both families are separating, for any point set,
hence they contain
separating subfamilies of cardinality at most $n-1$. This implies by
Observation \ref{mini} an upper bound of $n-1$ in each case.
In fact, we can give a separating subfamily of size $n-1$ directly for the hyperplanes. Assume without loss of generality that all points have different $x$-coordinates,
(otherwise we slightly rotate the coordinate system)
and let $p_1, p_2, \ldots , p_n$ be the points, ordered according to their
$x$-coordinates.
Then let $H_1, H_2, \ldots , H_{n-1}$ be halfplanes such that $H_i$ contains
exactly $p_1, p_2, \ldots , p_i$ of the points. It is clearly a separating family.

Now let $P$ be a set of $n$ collinear points, $p_1, p_2, \ldots , p_n$, and let
$\mathcal{H}(P)$ be a separating family of halfplanes.
For any $i$, $1\le i<n$, there is a halfplane $H_i\in \mathcal{H}(P)$ which separates
$p_i$ and $p_{i+1}$, and these halfplanes are all different. Therefore,
$\mathcal{H}(P)$ contains at least $n-1$ halfplanes, consequently, $s(n,\mathcal{H})=n-1$.

Now we show the lower bound for $s'(n,\cal{H})$.
Let $P$ be a set of $n$ points, $p_1, p_2, \ldots , p_n$, on a circle, ordered clockwise,
and let
$\mathcal{H}(P)$ be a separating family of halfplanes.
For any $i$, $1\le i<n$, there is a halfplane $H_i\in \mathcal{H}(P)$ which separates
$p_i$ and $p_{i+1}$, and there is a  halfplane $H_n\in \mathcal{H}(P)$ which separates
$p_n$ and $p_{1}$. These are $n$ halfplanes, and we counted a halfplane at most twice.
Therefore,
$\mathcal{H}(P)$ contains at least ${\lceil n/2\rceil}$ halfplanes.

Finally, we show that the upper bound holds for  $s'(n,\cal{H})$.
Suppose that $P$ is a set of $n$ points in general position.
We obtain a separating family $\cal S$ of ${\lceil n/2\rceil}$ halfplanes by the following procedure.
We can assume without loss of generality that $n$ is even.

\bigskip

{\em {\sc Halfplane-Separate($P$)}

\smallskip

\noindent{\sc Step $0$.} Let $\ell$ be a line which has exactly $n/2$ of the points on both sides.
Let $Q_0\subset P$ and $R_0\subset P$
denote  the points on the two sides of $\ell$. Let $H_0$ be a halfplane whose
bounding line is $\ell$. Let ${\cal S}_0=\{ H_0\}$, $i=1$.

\noindent{\sc Step $i$.} Take the convex hull of $Q_i\cup R_i$. It has two edges that cross $\ell$, let $e=q_ir_i$
be one of them, $q_i\in Q_i$, $r_i\in R_i$. Take a halfplane $H_i$ that
separates $q_i$ and $r_i$
from the rest of $Q_i$ and $R_i$.

If $i<n/2-1$, then let ${\cal S}_{i+1}={\cal S}_{i}\cup \{ H_i\}$,
let $Q_{i+1}=Q_i\setminus \{ q_i\}$, $R_{i+1}=R_i\setminus \{ r_i\}$.
Increase $i$ by one, and 
go to {\sc Step $i$.}

Otherwise (for  $i=n/2-1$), let ${\cal S}={\cal S}_{i}\cup \{ H_i\}$, and {\sc Stop}.}

\bigskip

Clearly, $\cal S$ is a set of $n/2$ halfplanes. We claim that
it separates $P$. Let $p, p'\in P$. If $p=q_i$ and $p'=r_j$ for some $i, j$, then $H_0$ separates them.
If $p=q_i$ and $p'=q_j$, $i<j$, (or if  $p=r_i$ and $p'=r_j$, $i<j$),
then $H_i$ separates them.


In the case of discs, the proofs are very similar.
To show that  $s(n,\mathcal{D}), s'(n,\mathcal{D})\ge {\lceil n/2\rceil}$,
let $P$ be a set of $n$ points on a circle.
We can argue exactly as in the case of halfplanes, any disc can separate
at most two consecutive pairs, therefore we need at least
${\lceil n/2\rceil}$ discs.

To prove that ${\lceil n/2\rceil}$ discs are always enough, we can use
a procedure very similar to
{\sc Halfplane-Separate($P$)}. Observe, that in the case of discs, it
  works for any point set, we do not have to assume that the points are in
  general position. Therefore, we have  $s(n,\mathcal{D})=s'(n,\mathcal{D})={\lceil n/2\rceil}$.

\section{Convex sets; Proof of Theorem \ref{fo}}\label{3}

Proof of part (i). The proof
of the lower bound is again very similar to the previous proofs.
Let $P$ be a set of $n$ points, $p_1, p_2, \ldots , p_n$, on a line in this
order,
and let
$\mathcal{A}(P)$ be a separating family of convex sets.
For any $i$, $1\le i<n$, there is a set $A_i\in \mathcal{A}(P)$ which separates
$p_i$ and $p_{i+1}$, and there is a set $A_n\in \mathcal{A}(P)$ which separates
$p_n$ and $p_{1}$. These are $n$ sets, and we counted a set at most twice.
Therefore,
$\mathcal{A}(P)$ contains at least ${\lceil n/2\rceil}$ sets,
so $s(n,\mathcal{D})\ge {\lceil n/2\rceil}$.
On the other hand, since discs are convex sets,
 $s(n,\mathcal{A})\le s(n,\mathcal{D})={\lceil n/2\rceil}$.
Therefore,  $s(n,\mathcal{A})={\lceil n/2\rceil}$.

Proof of part (ii).
Let $ES(k)$ denote the least integer such that among any $ES(n)$ points in
general position in the plane there are always $k$ in convex position.
In 1935, P. Erd\H os and G. Szekeres \cite{ES35}
showed that $ES(k)$ exists for every $k$,
and
$ES(k)\le {2k-4\choose k-2}+1$. The best known bounds for $ES(k)$ are

$$2^{k-2}+1\le ES(k)\le {2k-5\choose k-2}+1,$$
they were proved by
P. Erd\H os and G. Szekeres \cite{ES60}, and by
T\'oth and Valtr \cite{TV05}, respectively.

It is easy to see that both the lower and upper bound holds for $n \le 16$, hence we can assume that $n>16$.

First we prove the lower bound for  $s'(n,\mathcal{A})$.
Assume without loss of generality that $n$ is even.
Using the construction of Erd\H os and Szekeres \cite{ES60}, or a subset of
it,
we can obtain a point set $P_n$ 
of size $n/2$, in general position, such that it does not contain more than
$2\log n$ points in convex position. Take a line $\ell$ which is not parallel
to any line determined by the points of $P_n$.
Substitute each pont $p$ of $P_n$ by two points, $p'$ and $p''$, both very
close to $p$ such that
the line $p'p''$ is parallel to $\ell$. Points $p'$ and $p''$ are called {\em twins}, and
$p$ is their {\em parent}.
Let $Q_n$ be the resulting set of $n$ points, which is clearly in general position.

Suppose that $\cal S$ separates $Q_n$.
Clearly, for each pair of twins $(p', p'')$ in $Q_n$, there is a set $S\in {\cal S}$
which separates
them, that is, is contains exactly one of $p'$ and $p''$. Assign such a set $S(p',p'')$ to each pair
$(p',p'')$.

This way we found $n/2$ members of $\cal S$. Now we estimate how many times we could find the same
set. Suppose e. g. that $S(p_1',p_1'')=S(p_2',p_2'')=\cdots =S(p_k',p_k'')$.
Then, since $S$ is convex, and twins are very close to each other and to their parents,
points $p_1, p_2, \ldots , p_k$ are in convex position.
Therefore, by the assumption, $k\le 2\log n$. So, the number of different sets assigned to
the twins is at least $n/(2\log n)$.

\smallskip

Now we prove the upper bound. Again, assume that $n>16$.

%
%


By \cite{TV05}, any set of $m$ points in
general position contains
$\log m/2$ in convex position.
Let $P$ be a set of $n$ points in general position. We select a
separating system $\cal S$ of convex sets such that its cardinality is at most
$20n\log\log n/\log n$
by the following procedure.

\bigskip

{\em {\sc Convex-Separate($P$)}

\smallskip

Let $P_1=P$, ${\cal S}_1=\emptyset$, $i=1$.

\noindent{\sc Step $i$}. Let $Q_i\subset P_i$ be a subset of $k=\lfloor\frac{ \log n}{4}\rfloor $
points in convex position. Then there is a family $\mathcal{A}_i$ of cardinality $\lceil \log k \rceil$ which separates $Q_i$.

Let $S_i$ be the convex hull of all points of $Q_i$, and
let $S'_i$ be slightly shrinked copy of $S_i$ (one which contains all the points of $P_i\cap S_i$, except for the points of $Q_i$).

Add $S_i$, $S'_i$, and $\mathcal{A}_i$ to $\cal S$.

Let $P_{i+1}= P_i \setminus Q_i$. If $|P_{i+1}|> \sqrt{n}$, then     
increase $i$ by one and
go to {\sc Step $i$.}

Otherwise, go to {\sc Final Step.}

\medskip

\noindent{\sc Final Step.} For each point $p\in P$, add $S(p)$,
a very small disc
with center $p$, to $\cal S$.

\medskip

{\sc Stop}.}

\bigskip

When we execute {\sc Step $i$},
$P_i$ contains more than $\sqrt{n}$ points,
hence, by \cite{TV05},  we can select $k=\lfloor\frac{\log
  n}{4}\rfloor>\frac{\log n}{5}$ points among them in convex position. In each step,
except for the final one,
we delete $k$ points, hence we repeat {\sc Step $i$} at most  $5n/\log n$ times.
Each time
we select at most $2\log\log n+2$ sets, and in the Final Step we select at most $\sqrt{n}$
sets to $\cal S$.
So, eventually, we have $|{\cal S}|\le 10n\log\log n/\log n+10n/\log
n+\sqrt{n}\le 20n\log\log n/\log n$.
We claim that $\cal S$ separates $P$. We have to show that
any two elements of $P$ can be separated by some member of
$\cal S$.
If $p, p'\in Q_i$ for some $i$, then $\mathcal{A}_i$ separates them.

Now suppose that $p\in Q_i$ and $p'\not\in Q_i$ for some $i$.
If $p'$ is in the convex hull of $Q_i$, then $S_0$ from
{\sc Step $i$} separates $p$ and $p'$,
if $p'$ is not in the convex hull of $Q_i$, then $S'_0$ from
{\sc Step $i$} separates $p$ and $p'$.

Finally, suppose that there is no $i$ such that $p\in Q_i$.
Then we selected set $S(p)$ in the {\sc Final Step}, and it separates
$p$ and $p'$.

This concludes the proof of Theorem \ref{fo}.

\bigskip

\noindent {\bf Remark.} Let
$\mathcal{A}$ be a family of connected planar sets, $\gamma$ a Jordan curve, and $k$ a constant.
Suppose that each $A\in{\mathcal{A}}$ is bounded by a closed Jordan curve, and
intersects $\gamma$ in at most $k$ intervals. 
%
%
Then, it is not hard to see \cite{P12} that 
 $s(n,\mathcal{A})\ge (n-1)/2k$. 
We believe that there is a linear bound for other ``simple'' families, in
particular, 
for families of finite Vapnik-Chervonenkis dimension.

\bigskip

\noindent {\bf Conjecture.} Suppose that $\mathcal{A}$ is a family of planar
sets, whose 
Vapnik-Chervonenkis dimension is finite. Then there is a $c=c(\mathcal{A})>0$
such that  $s(n,\mathcal{A})>cn$ for every $n$.

\bigskip

Note that this conjecture has nothing to do with geometry, it is a purely
combinatorial statement. On the other hand, it might be 
easier
to verify the conjecture if assume that the sets in 
 $\mathcal{A}$ 
are {\em connected}, or 
we add some geometric condition.

\end{document}